\documentclass[a4paper,12pt]{amsart}

\usepackage{amsmath}
\usepackage{amssymb}
\usepackage{mathrsfs}
\usepackage{ifthen}
\usepackage{graphicx}
\usepackage[T1]{fontenc} 

\def\switchlinenumbers{\@ifstar
    {\let\makeLineNumberOdd\makeLineNumberRight
     \let\makeLineNumberEven\makeLineNumberLeft}%
    {\let\makeLineNumberOdd\makeLineNumberLeft
     \let\makeLineNumberEven\makeLineNumberRight}%
    }

\def\setmakelinenumbers#1{\@ifstar
  {\let\makeLineNumberRunning#1%
   \let\makeLineNumberOdd#1%
   \let\makeLineNumberEven#1}%
  {\ifx\c@linenumber\c@runninglinenumber
      \let\makeLineNumberRunning#1%
   \else
      \let\makeLineNumberOdd#1%
      \let\makeLineNumberEven#1%
   \fi}%
  }

\nonstopmode \numberwithin{equation}{section}
\newtheorem*{theorem*}{Theorem}

\newtheorem{thm}{Theorem}[section]
\newtheorem{cor}[equation]{Corollary}
\newtheorem{lem}{Lemma}[section]

\theoremstyle{definition}

\newtheorem{prob}[equation]{Problem}

\newcounter{minutes}
\divide\time by 60
\newcounter{hours}
\multiply\time by 60
\addtocounter{minutes}{-\time}

\newcounter {own}
\def\theown {\thesection       .\arabic{own}}

\newenvironment{pf}[1][]{%
 \vskip 3mm
 \noindent
 \ifthenelse{\equal{#1}{}}%
  {{\slshape Proof. }}%
  {{\slshape #1.} }%
 }%
{\qed\bigskip}

\newcounter{alphabet}

\def\be{\begin{equation}}
\def\ee{\end{equation}}

\newcommand{\bee}{\begin{enumerate}}
\newcommand{\eee}{\end{enumerate}}

\newcommand{\blem}{\begin{lem}}
\newcommand{\elem}{\end{lem}}
\newcommand{\bthm}{\begin{thm}}
\newcommand{\ethm}{\end{thm}}
\newcommand{\bcor}{\begin{cor}}
\newcommand{\ecor}{\end{cor}}
\newcommand{\beg}{\begin{examp}}
\newcommand{\eeg}{\end{examp}}
\newcommand{\begs}{\begin{examples}}
\newcommand{\eegs}{\end{examples}}
\newcommand{\bdefe}{\begin{defin}}
\newcommand{\edefe}{\end{defin}}
\newcommand{\bprob}{\begin{prob}}
\newcommand{\eprob}{\end{prob}}
\newcommand{\bei}{\begin{itemize}}
\newcommand{\eei}{\end{itemize}}

\newcommand{\norm}[1]{\left\lVert#1\right\rVert}

\newcommand{\innpdct}[1]{\left\langle#1\right\rangle}

\begin{document}

\title{Composition-Differentiation Operator on Weighted Bergman Spaces}

\author{Vasudevarao Allu}
\address{Vasudevarao Allu,
School of Basic Sciences,
Indian Institute of Technology Bhubaneswar,
Bhubaneswar-752050, Odisha, India.}
\email{avrao@iitbbs.ac.in}

\author{Himadri Halder}
\address{Himadri Halder,
School of Basic Sciences,
Indian Institute of Technology Bhubaneswar,
Bhubaneswar-752050, Odisha, India.}
\email{himadrihalder119@gmail.com}

\author{Subhadip Pal}
\address{Subhadip Pal,
	School of Basic Sciences,
	Indian Institute of Technology Bhubaneswar,
	Bhubaneswar-752050, Odisha, India.}
\email{subhadippal33@gmail.com}

\subjclass[{AMS} Subject Classification:]{Primary 47B38, 47B33, 30H20; Secondary 47A05, 47B15}
\keywords{Composition Operator, Differentiation Operator, Weighted Bergman Space, Complex Symmetry, Derivative Hardy space}

\def\thefootnote{}
\footnotetext{ {\tiny File:~\jobname.tex,
printed: \number\year-\number\month-\number\day,
          \thehours.\ifnum\theminutes<10{0}\fi\theminutes }
} \makeatletter\def\thefootnote{\@arabic\c@footnote}\makeatother

\begin{abstract}
	In this paper, we study the complex symmetry of weighted composition-differentiation operator $D_{n, \psi, \phi}$ on weighted Bergman spaces $\mathcal{A}^2_{\alpha}$ with respect to the conjugation $C_{\mu, \eta}$ for $\mu, \eta \in \{z\in \mathbb{C}:|z|=1\}$. We obtain explicit conditions for which the operator $D_{n, \psi, \phi}$ is Hermitian and normal. We also characterize the complex symmetric weighted composition-differentiation operator for derivative Hardy spaces.
\end{abstract}

\maketitle
\pagestyle{myheadings}
\markboth{Vasudevarao Allu, Himadri Halder, and Subhadip Pal}{Composition-Differentiation Operator on Weighted Bergman Spaces}

\section{Introduction}
	The study of the composition operators has been initiated with two goals in
mind, namely, to have concrete examples of bounded operators on Hilbert spaces or
Banach spaces of functions, and to encounter the “Invariant Subspace Problem” of functional analysis from a different angle.	Banach and Stone (see \cite{Singh-book-1993}) long back ago have proved that linear, surjective isometries between
spaces of continuous functions on compact Hausdorff spaces are of type of a weighted composition operator. There have been several interesting generalizations and extensions of the Banach-Stone theorem by several auhors in many different contexts. Among the existing works, the isometries of function spaces, particularly on Hardy spaces $\mathcal{H}^p$ are much more interesting (see \cite{Forelli-1964,Leeuw-PAMS-1960}). Forelli \cite{Forelli-1964} and Leeuw \cite{Leeuw-PAMS-1960} have shown that, for $p\geq 1$, the only surjective isometries of Hardy spaces $\mathcal{H}^p$ are precisely type of weighted composition operator, except the case $p=2$. Banach-Stone theorem not only motivates and inspires to study isometries on function spaces, but also sets the stage for the study of weighted composition operators.

\vspace{2mm}
Let $\mathcal{H}(\mathbb{D})$ be the space of all analytic functions in the unit disk $\mathbb{D}:=\{z\in \mathbb{C}: |z|<1\}$. The Hardy space $\mathcal{H}^2$ is the Hilbert space consisting of functions $f\in \mathcal{H}(\mathbb{D})$ such that
\begin{equation*}
	\norm{f}^2_{\mathcal{H}^2}=\sup_{0<r<1}\frac{1}{2\pi} \int_{0}^{2\pi}|f(re^{i\theta})|^2\,d\theta<\infty.
\end{equation*}
For $-1<\alpha<\infty$, the weighted Bergman space $\mathcal{A}^2_{\alpha}$ is the set of all analytic functions on $\mathbb{D}$ such that 
\begin{equation*}
	\norm{f}^2_{\mathcal{A}^2_{\alpha}}=\int_{\mathbb{D}}|f(z)|^2(\alpha+1)(1-|z|^2)^{\alpha} \, dA(z)<\infty,
\end{equation*}
where $dA(z)$ is the normalized area measure on $\mathbb{D}$. When $\alpha=0$, the space is called the classical Bergman space, and it is denoted by $\mathcal{A}^2$ (see \cite{Hedenmalm-book-2000}). It is well-known that $\mathcal{A}^2_{\alpha}$ is a reproducing kernel Hilbert space with the reproducing kernel function at each $w\in \mathbb{D}$, 
\begin{equation*}
	K_w(z)=\frac{1}{(1-\overline{w}z)^{\alpha+2}},\,\,\, z\in \mathbb{D}.
\end{equation*}
For every non-negative integer $n$ and $w\in \mathbb{D}$, let $K^{[n]}_w(z)$ be the reproducing kernel for point-evaluation of the $n$-th derivative, where
\begin{equation*}
	K^{[n]}_w(z)=\frac{(\alpha+2)(\alpha+3)\cdots(\alpha+n+1)z^n}{(1-\overline{w}z)^{n+\alpha+2}}, \,\, z\in \mathbb{D}.
\end{equation*}
It can be shown that $\innpdct{f(z), K^{[n]}_w(z)}=f^{(n)}(w)$, for all $f \in \mathcal{A}^2_{\alpha}$. For any non-negative integer $n$, let 
\begin{equation*}
	\gamma_n=\sqrt{\frac{\Gamma(n+2+\alpha)}{n!\Gamma(2+\alpha)}}z^n, \,\,\,z\in \mathbb{D},
\end{equation*}
where $\Gamma(\cdot)$ denotes the standard Gamma function. Here the set $\{\gamma_n: n\geq 0\}$ forms an orthonormal basis for $\mathcal{A}^2_{\alpha}$ (see \cite{Hedenmalm-book-2000}).\\[3mm]

 Let $\psi \in \mathcal{H}(\mathbb{D})$ and $\phi$ be an analytic self-map of $\mathbb{D}$. A weighted composition operator on $\mathcal{H}(\mathbb{D})$ is a linear operator $\psi C_{\phi}$ given by 
\begin{equation*}
	(\psi C_{\phi} f)(z)= \psi(z)f(\phi(z)), \, f \in \mathcal{H}(\mathbb{D}).
\end{equation*} 
If $\psi(z)=1$, $\psi C_{\phi}$ reduces to the usual composition operator $C_{\phi}$. It is worth to note that every composition operator is bounded on $\mathcal{H}^2$ (see \cite{cowen-1995}). Composition operators on various analytic function spaces in $\mathbb{D}$ have been extensively studied by several authors. For more detailed exposition on composition operator, we refer to the books by Cowen and MacCluer \cite{cowen-1995} and Hedenmalm {\it et al.} \cite{Hedenmalm-book-2000} 

\vspace{2mm}
In the context of functions in $\mathcal{H}(\mathbb{D})$, it is also worth to consider the operators on $\mathcal{H}(\mathbb{D})$ defined in terms of differentiation. For $n\in \mathbb{N}$, the differential operator of order $n$ on $\mathcal{H}(\mathbb{D})$ is defined by $D^{(n)}(f)=f^{(n)}$. It is easy to see that the differentiation operator of order $n$ is not bounded on $\mathcal{A}^2 _{\alpha}$. However, for many analytic self-maps $\phi$ of $\mathbb{D}$, the operator $C_{\phi}D^{(n)}$ is given by 
\begin{equation*}
	f(z) \rightarrow f^n (\phi(z))
\end{equation*}
is bounded on $\mathcal{A}^2 _{\alpha}$. Hibschweiler and Portnoy \cite{hibschweiler-2005} and Ohno \cite{ohno-2006} have initially considered these operators and afterwards these operators have been noticed by many researchers (see \cite{fatehi-PAMS-2020, fatehi-CAOT-2021, Gu-2018, Han-wang-CAOT-2021, Stevic-2009}). Boundedness and compactness of $C_{\phi}D^{(1)}$ on Hardy space and $C_{\phi}D^{(n)}$ on weighted Bergman spaces respectively have been studied by Ohno \cite{ohno-2006} and Stevi\'{c} \cite{Stevic-2009}. We denote the bounded operator $C_{\phi}D^{(n)}$ by $D_{n, \phi}$. For $m \in \mathbb{N} \cup \{0\}$ and $\psi \in \mathcal{H}(\mathbb{D})$, the weighted composition differentiation operator $D _{m, \psi,\phi}$ is defined by 
\begin{equation*}
	(D_{m, \psi,\phi} f)(z):= \psi(z) f^m (\phi(z)), \,\,\, f \in \mathcal{H}(\mathbb{D}).
\end{equation*}
For $m=0$, the operator $D_{m, \psi,\phi}$ reduces to $\psi C_{\phi}$. When $m=0$ and $\psi(z)=1$, $D_{m, \psi,\phi}$ reduces to the composition operator $C_{\phi}$. For $m=1$ and $\psi(z)=1$, we have $D_{m, \psi,\phi}=C_{\phi}D$. When $m=1$ and $\psi(z)=\phi '(z)$, $D _{m, \psi,\phi}$ becomes $D C_{\phi}$.

\vspace{2mm}
A \textit{conjugation} on a separable Hilbert space $\mathcal{H}$ is an antilinear map $\mathcal{C}: \mathcal{H} \rightarrow \mathcal{H}$ which satisfies $\innpdct{\mathcal{C}x, \mathcal{C}y}=\innpdct{y,x}$ for all $x,y \in \mathcal{H}$ and $\mathcal{C}^2=I$. We say that a bounded linear operator $T$ on $\mathcal{H}$ is \textit{complex symmetric} if there exists a conjugation $\mathcal{C}$ on $\mathcal{H}$ such that $T=\mathcal{C} T^{*} \mathcal{C}$. 
 \vspace{2mm}
 
 The concept of complex symmetric operators arises as a natural generalization of complex symmetric matrices. Since the last two decades, the theory complex symmetric operators has been extensively studied by many mathematicians. The development has been initiated by Garcia and Putinar \cite{Garcia-TAMS-2006, Garcia-TAMS-2007}, Garcia and Poore \cite{Garcia-JFA-2013} and Garcia and Wogen \cite{Garcia-JFA-2009, Garcia-TAMS-2010} in this direction. In 2014, Jung {\it et al.} \cite{Jung-JFA-2014} obtained some characterizations for complex symmetric weighted composition-differentiation operator on $H^2(\mathbb{D})$. In addition, the authors have studied the spectral properties of these operators in \cite{Jung-JFA-2014}. In 2014, Garcia and Hammond \cite{Garcia-OTAA-2014} investigated the complex symmetric weighted composition operator for weighted Hardy spaces. For more intriguing aspects on complex symmetric operators, we refer to the work of Fatehi \cite{fatehi-CVEE-2019}, Lim and Khoi \cite{Lim-Khoi-JMAA-2018}, Wang and Xao \cite{Wang-IJM-2016} and Wang and Han \cite{Wang-JMAA-2019}.\\[3mm]
  
 Let $\mathcal{S}^2(\mathbb{D})$ be the space of analytic functions whose derivatives are in Hardy space $\mathcal{H}^2$, and it is defined by 
 $$\mathcal{S}^2(\mathbb{D})=\left\{f\in \mathcal{H}(\mathbb{D}): \norm{f}^2_{\mathcal{S}^2}=|f(0)|^2+\norm{f'}^2_{\mathcal{H}^2}=|f(0)|^2+\sum_{m=1}^{\infty}m^2|f_m|^2<\infty\right\}. $$
 In 2017, Gu and Luo \cite{Gu-2018} introduced an equivalent norm on $\mathcal{S}^2(\mathbb{D})$ and defined another derivative Hardy space $\mathcal{S}^2_{1}(\mathbb{D})$ by
 \begin{align*}
 	\mathcal{S}^2_{1}(\mathbb{D})& =\left\{f\in \mathcal{H}(\mathbb{D}): \norm{f}^2_{\mathcal{S}^2_{1}}=\norm{f}^2_{\mathcal{H}^2}+\frac{3}{2}\norm{f'}^2_{\mathcal{A}^2}+\frac{1}{2}\norm{f'}^2_{\mathcal{H}^2}<\infty\right\} \\ \\[1mm] & = 
 	\left\{f\in \mathcal{H}(\mathbb{D}): \norm{f}^2_{\mathcal{S}^2_{1}}= \sum_{m=0}^{\infty}\frac{(m+1)(m+2)}{2}|f_m|^2< \infty\right\}.
 \end{align*} 
 For $f, g \in \mathcal{S}^2_{1}(\mathbb{D})$, the inner product on $\mathcal{S}^2_{1}(\mathbb{D})$ is given by
 $$ \innpdct{f, g}_{\mathcal{S}^2_{1}}=\innpdct{f, g}_{\mathcal{H}^2}+\frac{3}{2}\innpdct{f', g'}_{\mathcal{A}^2}+\frac{1}{2}\innpdct{f', g'}_{\mathcal{H}^2}.$$
 For any $w\in \mathbb{D}$ define 
 $$ K_w(z)=\sum_{m=0}^{\infty}\delta_{m}(\overline{w}z)^m, \,\, z\in \mathbb{D},$$
 where $$\delta_{m}=\frac{2}{(m+1)(m+2)} \quad \,\, \mbox{for}\,\, m \in \mathbb{N}\cup\{0\}.$$
 It is easy to see that
 \begin{equation*}
 	\innpdct{f(z), K_w(z)}_{\mathcal{S}^2_{1}}=f(w), \quad z\in \mathbb{D}.
 \end{equation*}
 Therefore, $K_w$ is the reproducing kernel function for $\mathcal{S}^2_{1}(\mathbb{D})$.
 For every $w\in \mathbb{D}$ and $n\in \mathbb{N}$, we define the reproducing kernel for point evaluation of the $n$-th derivative as $$K^{[n]}_w(z)=\sum_{m=n}^{\infty}\frac{m!}{(m-n)!}\delta_{m}(\overline{w})^{m-n}z^m$$
 such that $\innpdct{f, K^{[n]}_w}=f^{(n)}(w)$ for all $z\in \mathbb{D}$. Moreover, we observe that for any $f \in \mathcal{S}^2_{1}(\mathbb{D})$,
 \begin{equation*}
 	\innpdct{f, D^*_{n, \psi, \phi}K_w(z)}_{\mathcal{S}^2_{1}(\mathbb{D})}=\innpdct{D_{n, \psi, \phi}f, K_w(z)}_{\mathcal{S}^2_{1}(\mathbb{D})}=\psi(w)f^{(n)}(\phi(w))=\innpdct{f, \overline{\psi(w)}K^{[n]}_{\phi(w)}}_{\mathcal{S}^2_{1}(\mathbb{D})}.
 \end{equation*}
 Therefore, $D^*_{n, \psi, \phi}K_w(z)=\overline{\psi(w)}K^{[n]}_{\phi(w)}$.\\[2mm]
 
  The main aim of this paper is to study complex symmetric weighted composition-differentiation operator on weighted Bergman space $\mathcal{A}^2_{\alpha}$ and on the derivative Hardy space $\mathcal{S}^2_{1}(\mathbb{D})$ with respect to the conjugation $C_{\mu, \eta}$. For $\mu, \eta$ on the unit circle $\{z\in \mathbb{C}:|z|=1\}$, $C_{\mu, \eta}$ is defined by
 \begin{equation*}
 	C_{\mu, \eta}f(z)=\mu \overline{f(\overline{\eta z})}
 \end{equation*}
 for all analytic function $f$. The organization of this paper is as follows. In Theorem \ref{thm-2.1}, we characterize the functions $\psi$ and $\phi$ such that $D_{n, \psi, \phi}$ is complex symmetric on $A^2_{\alpha}$ with conjugation $C_{\mu, \eta}$. When $\phi$ is an automorphism on $\mathbb{D}$, Theorem \ref{thm-2.2} gives explicit forms of $\phi$ such that $D_{n, \psi, \phi}$ is bounded and $C_{\mu, \eta}$-symmetric. In Theorem \ref{thm-2.3}, we obtain a sufficient condition for $D_{n, \psi, \phi}$ to be normal whereas, Theorem \ref{Vasu-Him-Sub-P2-e-2.4} concerns a necessary and sufficient condition for $D_{n, \psi,\phi}$ to be Hermitian. In Theorem \ref{thm-2.5}, we characterize the functions $\psi$ and $\phi$ such that $D_{n, \psi, \phi}$ is complex symmetric on $A^2_{\alpha}$ with conjugation $C_{\mu, \eta}$. For $n=1$, Theorems \ref{thm-2.1}, \ref{thm-2.4} and \ref{thm-2.5} generalize the results which are recently proved by Liu {\it et al.} \cite{Liu-Ponnu-Linear-2022}.

\section{Main Results}
We study the weighted composition-differentiation operator on weighted Bergman space $\mathcal{A}^2_{\alpha}$. We first prove the following lemma to prove one of our main results. Throughout the article, we shall denote $p=(\alpha+2)(\alpha+3)\cdots (\alpha+n+1)$.

\begin{lem}\label{lm-1}
	Let $\phi$ be an analytic self-map of $\mathbb{D}$ and $\psi\in \mathcal{H}(\mathbb{D})$ such that $D_{n, \psi, \phi}$ is bounded on $\mathcal{A}^2_{\alpha}$. Then for any $w\in \mathbb{D}$, 
	$$ D^*_{n, \psi, \phi}K_w(z)=\overline{\psi(w)}K^{[n]}_{\phi(w)}.$$
\end{lem}

We  now give a complete characterization for the complex symmetric weighted composition-differentiation operator on $\mathcal{A}^2_{\alpha}$ with the conjugation $C_{\mu, \eta}$.
\begin{thm}\label{thm-2.1}
Let $\phi: \mathbb{D}\rightarrow \mathbb{D}$ be an analytic self-map on $\mathbb{D}$ and $\psi \in \mathcal{H}(\mathbb{D})$ with $\psi \neq 0$ such that $D_{n, \psi, \phi}$ is bounded on $\mathcal{A}^2_{\alpha}$. Then $D_{n, \psi, \phi}$ is complex symmetric on $\mathcal{A}^2_{\alpha}$ with conjugation $C_{\mu, \eta}$ if, and only if, 
\begin{equation*}
	\psi(z)=\frac{az^n}{(1-\eta b z)^{n+\alpha+2}} \,\,\,\,\, \mbox{and}\,\,\,\,\,\,\, \phi(z)=b+\frac{cz}{1-\eta b z}\,\,\, \mbox{for all}\,\,\, z\in \mathbb{D},
\end{equation*}
where $a, c \in \mathbb{C}$ and $b\in \mathbb{D}$.
\end{thm}

 In the following result, we obtain the condition on $\phi$ so that $\phi$ is an automorphism on $\mathbb{D}$ and the operator $D_{n, \psi, \phi}$ is complex symmetric with conjugation $C_{\mu, \eta}$ on $\mathcal{A}^2_{\alpha}$.
 
 \begin{thm}\label{thm-2.2}
 	Let $n \in \mathbb{N}$, $\phi$ be an automorphism on $\mathbb{D}$ and $\psi \in \mathcal{H}(\mathbb{D})$ be not identically zero such that $D_{n, \psi, \phi}$ is bounded and $C_{\mu, \eta}$-symmetric on $\mathcal{A}^2_{\alpha}$. Then one of the following statements hold:
 	\begin{enumerate}
 		\item[(a)] $\phi(z)=-\xi z$ with $|\xi|=1$ for some $\xi \in \mathbb{C}$.\\[2mm]
 		\item[(b)] $$\phi(z)=\frac{\overline{a_0}}{\eta a_0}\cdot \frac{a_0 -z}{1-\overline{a_0}z}\,\,\,\, \mbox{for some} \,\, a_0 \in \mathbb{D}\backslash \{0\}.$$
 	\end{enumerate}
 \end{thm}
\vspace{2mm}

We obtain a sufficient condition for the $C_{\mu, \eta}$-symmetric bounded operator $D_{n, \psi,\phi}$ to be normal. A bounded linear operator $T$ defined on a Hilbert space $\mathcal{H}$ is said to be normal if $TT^*=T^*T$ or, for any $x\in \mathcal{H}$, $\norm{Tx}=\norm{T^*x}$.

\begin{thm}\label{thm-2.3}
	Let $\phi$ be an analytic self-map of $\mathbb{D}$ with $\phi(0)=0$ and $\psi \in \mathcal{H}(\mathbb{D})$ be not identically zero such that for any $n\in \mathbb{N}$, the operator $D_{n, \psi, \phi}$ is bounded and complex symmetric with conjugation $C_{\mu, \eta}$. Then $D_{n, \psi, \phi}$ is normal.
\end{thm}
\vspace{2mm}

A bounded linear operator is called Hermitian if $T=T^*$. In the following result, we establish a necessary and sufficient condition for the bounded operator $D_{n, \psi, \phi}$ to be Hermitian on $\mathcal{A}^2_{\alpha}$, $n \in \mathbb{N}$. 
 \begin{thm}\label{thm-2.4}
	Let $\phi$ be an analytic self-map of $\mathbb{D}$ and $\psi$ be a nonzero analytic function on the unit disk $\mathbb{D}$ such that $D_{n, \psi, \phi}$ is bounded on $\mathcal{A}^2_{\alpha}$. Then $D_{n, \psi, \phi}$ is Hermitian if, and only if, 
	$$ \psi(z)=\frac{az^n}{(1-\overline{b}z)^{n+\alpha+2}} \,\,\, \mbox{and} \,\,\, \phi(z)=b+\frac{cz}{1-\overline{b}z}$$ 
	for some $a, c\in \mathbb{R}$ and $b \in \mathbb{D}$.
\end{thm}
\vspace{1mm}

Next, we study the weighted composition-differentiation operator on the derivative Hardy space $\mathcal{S}^2_{1}(\mathbb{D})$. We give a complete characterization of complex symmetric weighted composition-differentiation operator $D_{n, \psi, \phi}$ on $\mathcal{S}^2_{1}(\mathbb{D})$.

\begin{thm}\label{thm-2.5}
 Let $\phi$ be an analytic self-map on $\mathbb{D}$ and $\psi \in \mathcal{H}(\mathbb{D})$ such that $D_{n, \psi, \phi}$ is bounded on $\mathcal{S}^2_{1}(\mathbb{D})$. If $D_{n, \psi, \phi}$ is complex symmetric on $\mathcal{S}^2_{1}(\mathbb{D})$ with conjugation $C_{\mu, \eta}$, then there exist $a,b,c\in \mathbb{C}$ such that 
 \begin{equation}\label{Vasu-Him-Sub-P2-e-2.20}
 	\psi(z)=\frac{a}{n!\delta_{n}}F_1(z)\,\, \mbox{and}\,\, \phi(z)=b+c\frac{(n+3)}{(n+1)^2}\cdot\frac{F_2(z)}{F_1(z)}\,\,\,\, \mbox{for all}\,\, z\in \mathbb{D},
 \end{equation}
where $$ F_1(z)=\sum_{k=n}^{\infty}\frac{k!}{(k-n)!}\delta_{k}(\eta b)^{k-n}z^k\,\,\, \mbox{and}\,\,\, F_2(z)=\sum_{k=n+1}^{\infty}\frac{k!}{(k-n-1)!}\delta_{k}(\eta b)^{k-n-1}z^k.$$
Conversely, if $\psi$ and $\phi$ are defined by \eqref{Vasu-Him-Sub-P2-e-2.20} for some $a,b,c\in \mathbb{C}$, then $D_{n, \psi,\phi}$ is complex symmetric on $\mathcal{S}^2_{1}(\mathbb{D})$ with conjugation $C_{\mu, \eta}$ only if $b=0$ or $c=0$ or both are zero.
\end{thm}

\section{Proof of the Main Results}
\vspace{-2mm}

\begin{pf} [{\bf Proof of Lemma \ref{lm-1}}]
		Let $f\in \mathcal{A}^2_{\alpha}$. Then
	\begin{align*}
		\innpdct{f, D^*_{n, \psi, \phi}K_w(z)}_{\mathcal{A}^2_{\alpha}} &= \innpdct{D_{n, \psi, \phi}f, K_w}_{\mathcal{A}^2_{\alpha}}\\ &=\innpdct{\psi\cdot f^{(n)}(\phi)K_w}_{\mathcal{A}^2_{\alpha}}\\ &= \psi(w)f^{(n)}(\phi(w))=\innpdct{f, \overline{\psi(w)}K^{[n]}_{\phi(w)}}_{\mathcal{A}^2_{\alpha}}
	\end{align*}
	holds for any $f\in \mathcal{A}^2_{\alpha}$. Consequently, $ D^*_{n, \psi, \phi}K_w(z)=\overline{\psi(w)}K^{[n]}_{\phi(w)}$.
\end{pf}

\begin{pf} [{\bf Proof of Theorem \ref{thm-2.1}}]
		Let $D_{n, \psi, \phi}$ be complex symmetric with conjugation $C_{\mu, \eta}$. Then we have 
	\begin{equation}\label{Vasu-Him-Sub-P2-e-2.1}
		D_{n,\psi,\phi}C_{\mu, \eta} K_{w}(z)=C_{\mu, \eta}D^{*}_{n, \psi, \phi} K_{w}(z)
	\end{equation}
	for all $z, w \in \mathbb{D}$.
	A simple computation shows that 
	\begin{align}
		D_{n,\psi,\phi}C_{\mu, \eta} K_{w}(z) \label{V-H-P-e-2.2-b}& = D_{n,\psi,\phi}C_{\mu, \eta} \left(\frac{1}{(1-\overline{w}z)^{\alpha+2}}\right) \\\nonumber & = D_{n,\psi,\phi} \left(\frac{\mu}{(1-\eta wz)^{\alpha+2}}\right) \\\nonumber & = \frac{\mu p\psi(z)(\eta w)^n}{(1-\eta w\phi(z))^{n+\alpha+2}}.
	\end{align}
	Therefore, the right hand side of \eqref{Vasu-Him-Sub-P2-e-2.1} yields
	\begin{equation}\label{V-H-P-e-2.3-b}
		C_{\mu, \eta}D^{*}_{n, \psi, \phi} K_{w}(z)  = C_{\mu, \eta} \overline{\psi(w)} K^{[n]}_{\phi(w)} = \frac{\mu p\psi(w)(\eta z)^n}{(1-\eta \phi(w) z)^{n+\alpha +2}}.
	\end{equation}
	By making use of \eqref{V-H-P-e-2.2-b} and \eqref{V-H-P-e-2.3-b} in \eqref{Vasu-Him-Sub-P2-e-2.1}, we obtain
	\begin{equation}\label{Vasu-Him-Sub-P2-e-2.2}
		\frac{\mu p\psi(z)(\eta w)^n}{(1-\eta w\phi(z))^{n+\alpha+2}}= \frac{\mu p\psi(w)(\eta z)^n}{(1-\eta \phi(w) z)^{n+\alpha +2}}
	\end{equation}
	for all $z, w \in \mathbb{D}$. Letting $z=0$ in \eqref{Vasu-Him-Sub-P2-e-2.2}, we obtain
	\begin{equation}\label{V-H-P-e-2.3-a}
		\frac{\mu p\psi(0)(\eta w)^n}{(1-\eta w\phi(0))^{n+\alpha+2}}=0
	\end{equation}
	for any $z \in \mathbb{D}$. Since $\mu, \eta \in \{z\in \mathbb{C}:|z|=1\}$, it is easy to see that \eqref{V-H-P-e-2.3-a} is equivalent to $\psi(0)=0$.\\[2mm]
	
	Let $\psi(z)=z^mg(z)$, where $m\in \mathbb{N}$ and $g$ is analytic on $\mathbb{D}$ with $g(0)\neq 0$. Now our aim is to show that $m=n$. If $m>n$, from \eqref{Vasu-Him-Sub-P2-e-2.2} it follows that 
	\begin{equation}\label{Vasu-Him-Sub-P2-e-2.3}
		\frac{z^{m-n}g(z)}{(1-\eta w\phi(z))^{n+\alpha+2}}= \frac{w^{m-n}g(w)}{(1-\eta \phi(w) z)^{n+\alpha +2}}
	\end{equation}
	for any $z, w \in \mathbb{D}$. Putting $w=0$ in \eqref{Vasu-Him-Sub-P2-e-2.3} we obtain $g(0)=0$, which is a contradiction of the fact that $g(0)\neq 0$. If $m<n$, from \eqref{Vasu-Him-Sub-P2-e-2.2} we obtain 
	\begin{equation}\label{Vasu-Him-Sub-P2-e-2.4}
		\frac{w^{n-m}g(z)}{(1-\eta w\phi(z))^{n+\alpha+2}}= \frac{z^{n-m}g(w)}{(1-\eta \phi(w) z)^{n+\alpha +2}}.
	\end{equation}
	Setting $w=0$ in \eqref{Vasu-Him-Sub-P2-e-2.4} gives $g(0)=0$, which contradicts the assumption that $g(0)\neq 0$. Therefore, $m=n$ and hence \eqref{Vasu-Him-Sub-P2-e-2.2} reduces to 
	\begin{equation}\label{V-H-P-e-2.6-a}
		\frac{g(z)}{(1-\eta w\phi(z))^{n+\alpha+2}}= \frac{g(w)}{(1-\eta \phi(w) z)^{n+\alpha +2}}
	\end{equation}  
	for all $z, w \in \mathbb{D}$. By letting $w=0$ in \eqref{V-H-P-e-2.6-a} we obtain
	$$ g(z)=\frac{g(0)}{(1-\eta z \phi(0))^{n+\alpha+2}}.$$
	Therefore, we have
	\begin{align*}
		\psi(z)& =z^ng(z)\\ &= \frac{az^n}{(1-\eta bz)^{n+\alpha+2}}, \,\,\, \mbox{where}\,\,\, a=g(0) \,\, \mbox{and}\,\, b=\phi(0).
	\end{align*}
	By substituting  $\psi(z)$ in \eqref{Vasu-Him-Sub-P2-e-2.2}, we obtain 
	\begin{equation}\label{V-H-P-e-2.7-a}
		(1-\eta bz)^{n+\alpha+2}(1-\eta w\phi(z))^{n+\alpha+2}= (1-\eta bw)^{n+\alpha +2}(1-\eta \phi(w)z)^{n+\alpha+2}
	\end{equation}
	for all $z, w \in \mathbb{D}$. Differentiating both sides of \eqref{V-H-P-e-2.7-a} with respect to $w$, we obtain
	\begin{align}\label{Vasu-Him-Sub-P2-e-2.5}
		(n+\alpha+2) & (1-\eta bz)^{n+\alpha+2}  (1-\eta w\phi(z))^{n+\alpha+1}(-\eta \phi(z)) \\ \nonumber  & =(n+\alpha+2)(1-\eta bw)^{n+\alpha+1}(-\eta b)  \\ \nonumber & + (n+\alpha +2)(1-\eta bw)^{n+\alpha +2}(1-\eta \phi(w)z)^{n+\alpha+2}(-\eta z\phi'(w)).
	\end{align}
	Letting $w=0$ in \eqref{Vasu-Him-Sub-P2-e-2.5} we obtain
	\begin{equation*}
		\phi(z)=b+ \frac{cz}{1-\eta bz}
	\end{equation*}
	for all $z\in \mathbb{D}$, where $c=\phi'(0)$.\\[1mm]
	
	Conversely, let $a, c \in \mathbb{C}$ and $b\in \mathbb{D}$ be such that $$	\psi(z)=\frac{az^n}{(1-\eta b z)^{n+\alpha+2}} \,\,\,\,\, \mbox{and}\,\,\,\,\,\,\, \phi(z)=b+\frac{cz}{1-\eta b z}\,\,\, \mbox{for all}\,\,\, z\in \mathbb{D}. $$
	Then, we have
	\begin{align}
		D_{n,\psi,\phi}C_{\mu, \eta} K_{w}(z) \label{Vasu-Him-Sub-P2-e-2.6} & = \frac{\mu p\psi(z)(\eta w)^n}{(1-\eta w\phi(z))^{n+\alpha+2}}\\ \nonumber & = \frac{a\mu p (\eta zw)^n}{(1-\eta bz-\eta bw- \eta czw)^{n+\alpha +2}}.
	\end{align}
	On the other hand,
	\begin{align}
		C_{\mu, \eta}D^{*}_{n, \psi, \phi} K_{w}(z)\label{Vasu-Him-Sub-P2-e-2.7} & = C_{\mu, \eta} \overline{\psi(w)} K^{[n]}_{\phi(w)} \\ \nonumber & =  \frac{\mu p\psi(w)(\eta z)^n}{(1-\eta \phi(w) z)^{n+\alpha +2}}\\ \nonumber & = \frac{a\mu p(\eta zw)^n}{(1-\eta bw -\eta bz-\eta czw)^{n+\alpha+2}}.
	\end{align}
	In view of \eqref{Vasu-Him-Sub-P2-e-2.6} and \eqref{Vasu-Him-Sub-P2-e-2.7}, it follows that 
	$$D_{n,\psi,\phi}C_{\mu, \eta} K_{w}(z)=	C_{\mu, \eta}D^{*}_{n, \psi, \phi} K_{w}(z)\,\,\, \mbox{for all}\,\,\,z\in \mathbb{D}.$$ Since the span of the reproducing kernel functions is dense in $\mathcal{A}^2_{\alpha}$, the operator $D_{n, \psi, \phi}$ is complex symmetric with conjugation $C_{\mu, \eta}$. This completes the proof.
\end{pf}

\begin{pf}  [{\bf Proof of Theorem \ref{thm-2.2}}]
		Since $D_{n, \psi, \phi}$ is complex symmetric with conjugation $C_{\mu, \eta}$ on $\mathcal{A}^2_{\alpha}$, in view of Theorem \ref{thm-2.1} we have  
	$$ \phi(z)=b+ \frac{cz}{1-\eta bz},\,\,\, \mbox{where}\,\, \,b=\phi(0) \,\,\, \mbox{and}\,\,\, c=\phi'(0).$$
	Since $\phi$ is an automorphism on $\mathbb{D}$, there exist $a_0\in\mathbb{D}$ and $\xi \in \mathbb{C}$ with $|\xi|=1$ such that 
	\begin{equation}\label{Vasu-Him-Sub-P2-e-2.8}
		\phi(z)= b+ \frac{cz}{1-\eta bz}=\xi\cdot\frac{a_0-z}{1-\overline{a_0}z}
	\end{equation}
	for any $z\in \mathbb{D}$. It is easy to see that \eqref{Vasu-Him-Sub-P2-e-2.8} is equivalent to 
	\begin{equation}\label{Vasu-Him-Sub-P2-e-2.9}
		b-b\overline{a_0}z-\eta b^2z+\eta b^2\overline{a_0}z^2+cz-c\overline{a_0}z^2=\xi\overline{a_0}-\xi z -\eta b\xi a_0z+\xi\eta bz^2
	\end{equation}
	for any $z\in \mathbb{D}$. Comparing the constant terms of both sides of \eqref{Vasu-Him-Sub-P2-e-2.9} we obtain $b=\xi a_0$. Again, comparing the coefficients of $z$ and $z^2$, we have the following relations
	\begin{equation}\label{Vasu-Him-Sub-P2-e-2.10}
		c-b\overline{a_0}-\eta b^2=-\xi-\eta b\xi a_0 ,
	\end{equation}
	and 
	\begin{equation}\label{Vasu-Him-Sub-P2-e-2.11}
		\eta b^2\overline{a_0} - c\overline{a_0} = \xi \eta b.
	\end{equation}
	Now we consider the following two cases:
	
	{\it Case I:} If $a_0=0$, then $b=0$ and hence \eqref{Vasu-Him-Sub-P2-e-2.10} is equivalent to $c=-\xi$, which shows that 
	$$ \phi(z)=\xi z\,\,\, \mbox{with}\,\,\,\, |\xi|=1.$$ 
	
	{\it Case II:} If $a_0 \neq 0$, using $b=\xi a_0$ and \eqref{Vasu-Him-Sub-P2-e-2.11}, we obtain 
	$$ c=\frac{\xi^2\eta a_0}{\overline{a_0}}\cdot (|a_0|^2-1).$$
	By making use of the value of $c$ in \eqref{Vasu-Him-Sub-P2-e-2.10} we obtain $$ \xi=\frac{\overline{a_0}}{\eta a_0}.$$
	This completes the proof.
\end{pf}
\vspace{-2mm}

\begin{pf} [{\bf Proof of Theorem \ref{thm-2.3}}]
		Since $D_{n, \psi,\phi}$ is $C_{\mu, \eta}$-symmetric on $\mathcal{A}^2_{\alpha}$ and $\phi(0)=0$, in view of Theorem \ref{thm-2.1}, we have 
	$$\psi(z)=az^n \,\,\, \mbox{and}\,\,\, \phi(z)=cz,\,\,\,\, \mbox{where}\,\, a, c\in \mathbb{C}.$$
	Then for any non-negative integer $k$, we have
	\begin{align*}
		\norm{D_{n, \psi,\phi}\gamma_{k}}^2 &= \sum_{j=0}^{\infty}\left|\innpdct{D_{n, \psi,\phi}\gamma_{k}, \gamma_{j}}\right|^2 \\ & = \sum_{j=0}^{\infty}\left|\innpdct{\psi\cdot\gamma^{(n)}_{k}(\phi), \sqrt{\frac{\Gamma(j+2+\alpha)}{j!\Gamma(2+\alpha)}}z^j}\right|^2\\ & =
		\sum_{j=0}^{\infty}\left|\innpdct{\frac{k!ac^{k-n}}{(k-n)!}\sqrt{\frac{\Gamma(k+2+\alpha)}{k!\Gamma(2+\alpha)}}z^k, \sqrt{\frac{\Gamma(j+2+\alpha)}{j!\Gamma(2+\alpha)}}z^j }\right|^2
	\end{align*}
	and
	\begin{align*}
		\norm{D^*_{n, \psi, \phi}\gamma_{k}}^2 &= \sum_{j=0}^{\infty}\left|\innpdct{D^*_{n, \psi, \phi}\gamma_{k}, \gamma_{j}}\right|^2 \\ &= \sum_{j=0}^{\infty}\left|\innpdct{\gamma_{k}, D_{n, \psi, \phi}\gamma_{j}}\right|^2\\ &= \sum_{j=0}^{\infty}\left|\innpdct{\sqrt{\frac{\Gamma(k+2+\alpha)}{k!\Gamma(2+\alpha)}}z^k, \psi\cdot\gamma^{(n)}_{j}(\phi)}\right|^2 \\ &= \sum_{j=0}^{\infty}\left|\innpdct{\sqrt{\frac{\Gamma(k+2+\alpha)}{k!\Gamma(2+\alpha)}}z^k, \frac{j!ac^{j-n}}{(j-n)!}\sqrt{\frac{\Gamma(j+2+\alpha)}{j!\Gamma(2+\alpha)}}z^j}\right|^2.
	\end{align*}
	As the set $\{\gamma_n:n\geq 0\}$ forms an orthonormal basis for $\mathcal{A}^2_{\alpha}$ it follows that for any $k\in \mathbb{N}\cup\{0\}$, we have 
	\begin{equation}\label{Vasu-Him-Sub-P2-e-2.12-a}
		\norm{D_{n, \psi,\phi}\gamma_{k}}^2=\norm{D^*_{n, \psi, \phi}\gamma_{k}}^2=|ac^{k-n}|^2\left(\frac{k!}{(k-n)!}\right).
	\end{equation}
	Consequently, from \eqref{Vasu-Him-Sub-P2-e-2.12-a} it follows that $	\norm{D_{n, \psi,\phi}(f)}=\norm{D^*_{n, \psi, \phi}(f)}$ for all $f\in \mathcal{A}^2_{\alpha}$, and $D_{n, \psi, \phi}$ is a normal operator. This completes the proof.
\end{pf}
\begin{pf} [{\bf Proof of Theorem \ref{thm-2.4}}]
		Let $D_{n, \psi, \phi}$ be Hermitian. Then $D^*_{n, \psi, \phi}=D_{n, \psi, \phi}$, which is equivalent to 
	\begin{equation*}
		D^*_{n, \psi, \phi}K_{w}(z)=D_{n, \psi, \phi}K_{w}(z)\,\,\, \mbox{for all}\,\, z\in \mathbb{D}.
	\end{equation*}
	We observe that 
	\begin{equation}\label{V-H-P-e-2.16-a}
		D^*_{n, \psi, \phi}K_{w}(z)=\frac{pz^n\overline{\psi(w)}}{(1-z\overline{\phi(w)})^{n+\alpha+2}} \,\,\, \mbox{and}\,\,\, D_{n, \psi, \phi}K_{w}(z)=\frac{p(\overline{w})^n\psi(z)}{(1-\overline{w}\phi(z))^{n+\alpha+2}}.
	\end{equation}
	Therefore, from \eqref{V-H-P-e-2.16-a} we have 
	\begin{equation}\label{Vasu-Him-Sub-P2-e-2.12}
		\frac{pz^n\overline{\psi(w)}}{(1-z\overline{\phi(w)})^{n+\alpha+2}}=\frac{p(\overline{w})^n\psi(z)}{(1-\overline{w}\phi(z))^{n+\alpha+2}}. 
	\end{equation}
	for all $z, w \in \mathbb{D}$. Letting $w=0$ in \eqref{Vasu-Him-Sub-P2-e-2.12} we obtain $\psi(0)=0$. For $\psi \in \mathcal{A}^2_{\alpha}$, if $\psi(z)=z^mg(z)$, where $m \in \mathbb{N}$ and $g$ is analytic on $\mathbb{D}$ with $g(0)\neq 0$, then we prove that that $m=n$. 
	
	{\it Case I:} If $m>n$, then \eqref{Vasu-Him-Sub-P2-e-2.12} becomes 
	\begin{equation}\label{Vasu-Him-Sub-P2-e-2.13}
		\frac{z^{m-n}g(z)}{(1-\overline{w}\phi(z))^{n+\alpha+2}}=\frac{(\overline{w})^{m-n}\overline{g(w)}}{(1-z\overline{\phi(w)})^{n+\alpha+2}}
	\end{equation}
	for all $z, w \in \mathbb{D}$. Setting $w=0$ in \eqref{Vasu-Him-Sub-P2-e-2.13} it follows $g(0)=0$, which is a contradiction to the assumption that $g(0)\neq 0$.
	
	{\it Case II:} If $m<n$, then from \eqref{Vasu-Him-Sub-P2-e-2.12} we obtain 
	\begin{equation}\label{Vasu-Him-Sub-P2-e-2.14}
		\frac{z^{n-m}\overline{g(w)}}{(1-z\overline{\phi(w)})^{n+\alpha+2}}=\frac{\overline{w}^{n-m}g(z)}{(1-\overline{w}\phi(z))^{n+\alpha+2}}
	\end{equation}
	for all $z, w\in \mathbb{D}$. Putting $w=0$ in \eqref{Vasu-Him-Sub-P2-e-2.14} we obtain $g(0)=0$, which is a contradiction. Therefore, $m=n$. Now \eqref{Vasu-Him-Sub-P2-e-2.12} is equivalent to 
	\begin{equation}\label{V-H-P-e-2.20-a}
		\frac{\overline{g(w)}}{(1-z\overline{\phi(w)})^{n+\alpha+2}}=\frac{g(z)}{(1-\overline{w}\phi(z))^{n+\alpha+2}}
	\end{equation}
	for all $z, w\in \mathbb{D}$. Letting $w=0$ in \eqref{V-H-P-e-2.20-a}, it follows that 
	$$g(z)=\frac{\overline{g(0)}}{(1-z\overline{\phi(0)})^{n+\alpha+2}}.$$
	Therefore, 
	\begin{equation}\label{Vasu-Him-Sub-P2-e-2.15}
		\psi(z)=\frac{az^n}{(1-\overline{b}z)^{n+\alpha+2}}, \,\,\, \mbox{where}\,\,\, a=\overline{g(0)}\,\, \mbox{and}\,\, b=\phi(0).
	\end{equation}
	Substituting \eqref{Vasu-Him-Sub-P2-e-2.15} in \eqref{Vasu-Him-Sub-P2-e-2.12} we obtain
	\begin{equation}\label{Vasu-Him-Sub-P2-e-2.16}
		a(1-b\overline{w})^{n+\alpha+2}(1-z\overline{\phi(w)})^{n+\alpha+2}=\overline{a}(1-\overline{b}z)^{n+\alpha+2}(1-\overline{w}\phi(z))^{n+\alpha+2}
	\end{equation}
	for all $z, w\in \mathbb{D}$. If we set $w=0$ in \eqref{Vasu-Him-Sub-P2-e-2.16} we obtain that $a=\overline{a}$, and hence $a\in \mathbb{R}$. \\[2mm]
	
	Differentiating both sides of the equation \eqref{Vasu-Him-Sub-P2-e-2.16} with respect to $\overline{w}$ we obtain
	\begin{align}
		(n+\alpha+2)\label{V-H-P-e-2.23-a}& (1-\overline{b}z)^{n+\alpha+2}(1-\overline{w}\phi(z))^{n+\alpha+1}(-\phi(z)) \\\nonumber & = (n+\alpha+2)(1-b\overline{w})^{n+\alpha+1}(1-z\overline{\phi(w)})^{n+\alpha+2}(-b)\\ \nonumber & \quad \quad \quad + (n+\alpha+2)(1-b\overline{w})^{n+\alpha+2}(1-z\overline{\phi(w)})^{n+\alpha+1}(-z\phi'(w)) 
	\end{align}
	fol all $z\in \mathbb{D}$. For $w=0$, it is easy to see that \eqref{V-H-P-e-2.23-a} yields 
	\begin{equation}\label{Vasu-Him-Sub-P2-e-2.17}
		\phi(z)=b+\frac{cz}{1-\overline{b}z},\,\, z\in \mathbb{D},
	\end{equation}
	where $b=\phi(0)$ and $c=\phi'(0)$. In view of \eqref{Vasu-Him-Sub-P2-e-2.17}, we obtain $\phi'(0)=c$. Therefore, $c=\overline{\phi'(0)}=\overline{c}$ which implies that  $c \in \mathbb{R}$.\\[1mm]
	
	Conversely, assume that $\psi$ and $\phi$ be of the form given by \eqref{Vasu-Him-Sub-P2-e-2.15} and \eqref{Vasu-Him-Sub-P2-e-2.17}. Then we have 
	\begin{align}
		D_{n, \psi, \phi}K_w(z)\label{V-H-P-e-2.25-a} &=\frac{p(\overline{w})^n\psi(z)}{(1-\overline{w}\phi(z))^{n+\alpha+2}}\\\nonumber & = \frac{ap(z\overline{w})^n}{(1-\overline{b}z-b\overline{w}+\overline{w}|b|^2z-\overline{w}cz)^{n+\alpha+2}}
	\end{align} 
	and 
	\begin{align}
		D^*_{n, \psi,\phi}K_w(z)\label{V-H-P-e-2.26-a}&=\frac{apz^n\overline{\psi(w)}}{(1-z\overline{\phi(w)})^{n+\alpha+2}}\\ \nonumber& = \frac{ap(z\overline{w})^n}{(1-b\overline{w}-\overline{b}z+\overline{w}|b|^2z-\overline{w}cz)^{n+\alpha+2}}.
	\end{align}
	From \eqref{V-H-P-e-2.25-a} and \eqref{V-H-P-e-2.26-a}, it is clear that $D_{n, \psi, \phi}K_w(z)=	D^*_{n, \psi,\phi}K_w(z)$ for all $z\in \mathbb{D}$. Since the span of reproducing kernel functions is dense in $\mathcal{A}^2_{\alpha}$, it follows that $	D^*_{n, \psi,\phi}=	D_{n, \psi, \phi}$.
	This completes the proof.
\end{pf}
\begin{pf} [{\bf Proof of Theorem \ref{thm-2.5}}]
 Let $D_{n, \psi, \phi}$ be complex symmetric with conjugation $C_{\mu, \eta}$. Then we have the following:
\begin{equation}\label{V-H-P-e-2.28-a}
	D_{n, \psi, \phi}C_{\mu, \eta}K_w(z)=C_{\mu, \eta}	D^*_{n, \psi, \phi}K_w(z)
\end{equation}
for all $z, w\in \mathbb{D}$. Using \eqref{V-H-P-e-2.28-a}, we obtain
\begin{align*}
	D_{n, \psi, \phi}C_{\mu, \eta}K_w(z)& = 	D_{n, \psi, \phi}C_{\mu, \eta}\left(\sum_{k=0}^{\infty}\delta_{k}(\overline{w}z)^k\right)\\ & = 	D_{n, \psi, \phi}\left(\sum_{k=0}^{\infty}\mu\delta_{k}(\eta wz)^k\right)\\ &= \mu\psi(z) \sum_{k=n}^{\infty}\frac{k!}{(k-n)!}\delta_{k}(\eta w)^k\phi(z)^{k-n}.
\end{align*}
On the other hand, 
\begin{align*}
	C_{\mu, \eta}	D^*_{n, \psi, \phi}K_w(z)&= C_{\mu, \eta} \overline{\psi(w)}K^{[n]}_{\phi(w)}\\ &=  C_{\mu, \eta} \overline{\psi(w)} \sum_{k=n}^{\infty}\frac{k!}{(k-n)!}\delta_kz^{k}(\overline{\phi(w)})^{k-n} \\ &= \mu\psi(w)\sum_{k=n}^{\infty}\frac{k!}{(k-n)!}\delta_{k}(\eta z)^{k}\phi(w)^{k-n}.
\end{align*}
Therefore, we have
\begin{equation}\label{Vasu-Him-Sub-P2-e-2.21}
	\mu\psi(z) \sum_{k=n}^{\infty}\frac{k!}{(k-n)!}\delta_{k}(\eta w)^k\phi(z)^{k-n}=\mu\psi(w)\sum_{k=n}^{\infty}\frac{k!}{(k-n)!}\delta_{k}(\eta z)^{k}\phi(w)^{k-n}
\end{equation}
for all $z, w \in \mathbb{D}$. By substituting $w=0$ in \eqref{Vasu-Him-Sub-P2-e-2.21}, we obtain
$$ \mu \psi(0)\sum_{k=n}^{\infty}\frac{k!}{(k-n)!}\delta_{k}(\eta z)^k\phi(0)^{k-n}=0, \,\,\, z\in \mathbb{D}$$
which implies that $\psi(0)=0$. \\[2mm]

Let $\psi(z)=z^m\chi(z)$, where $m \in \mathbb{N}$ and $\chi(z)$ be an analytic function on $\mathbb{D}$ with $\chi(0)\neq 0$. By going the similar lines of argument as in the proof of Theorem \ref{thm-2.1} and Theorem \ref{thm-2.4}, it can be easily shown that $m=n$. Therefore, \eqref{Vasu-Him-Sub-P2-e-2.21} becomes
\begin{equation}\label{Vasu-Him-Sub-P2-e-2.22}
	\chi(z) \sum_{k=n}^{\infty}\frac{k!}{(k-n)!}\delta_{k}\eta ^k(w\phi(z))^{k-n}=\chi(w)\sum_{k=n}^{\infty}\frac{k!}{(k-n)!}\delta_{k}\eta^{k}(z\phi(w))^{k-n}
\end{equation} 
for all $z, w \in \mathbb{D}$. If we set $w=0$ in \eqref{Vasu-Him-Sub-P2-e-2.22}, we obtain
\begin{equation*}
	\chi(z)=\frac{a}{n!\delta_{n}}\sum_{k=n}^{\infty}\frac{k!}{(k-n)!}\delta_{k}(\eta bz)^{k-n}, \,\, z\in \mathbb{D},
\end{equation*}
where $a=\chi(0)$ and $b=\phi(0)$. Therefore, we have
\begin{equation}\label{Vasu-Him-Sub-e-2.23}
	\psi(z)=\frac{a}{n!\delta_{n}}\sum_{k=n}^{\infty}\frac{k!}{(k-n)!}\delta_{k}z^k(\eta b)^{k-n}, \,\, z\in \mathbb{D}.
\end{equation}
Substituting the expression of $\chi(z)$ in \eqref{Vasu-Him-Sub-P2-e-2.22}, we obtain
\begin{align*}
	\left(\sum_{k=n}^{\infty}\frac{k!}{(k-n)!}\delta_{k}(\eta bz)^{k-n}\right) & \left(\sum_{k=n}^{\infty}\frac{k!}{(k-n)!}\delta_{k}\eta ^k(w\phi(z))^{k-n}\right)\\ \\[2mm]\nonumber &= \left(\sum_{k=n}^{\infty}\frac{k!}{(k-n)!}\delta_{k}(\eta bw)^{k-n}\right)\left(\sum_{k=n}^{\infty}\frac{k!}{(k-n)!}\delta_{k}\eta^{k}(z\phi(w))^{k-n}\right).
\end{align*}
Differentiating both the sides of the above relation with respect to $w$, we obtain
\begin{align*}
	&\left(\sum_{k=n}^{\infty}\frac{k!}{(k-n)!}\delta_{k}(\eta bz)^{k-n}\right)\left(\sum_{k=n+1}^{\infty}\frac{k!}{(k-n)!}(k-n)\delta_{k}\eta^kw^{k-n-1}\phi(z)^{k-n}\right)\\ \\[3mm]&= \left(\sum_{k=n+1}^{\infty}\frac{k!}{(k-n)!}(k-n)\delta_{k}w^{k-n-1}(\eta b)^{k-n}\right)\left(\sum_{k=n}^{\infty}\frac{k!}{(k-n)!}\delta_{k}\eta^{k}(z\phi(w))^{k-n}\right)\\ \\[2mm]& \quad \quad \quad+ \left(\sum_{k=n}^{\infty}\frac{k!}{(k-n)!}\delta_{k}(\eta bw)^{k-n}\right)\left(\sum_{k=n+1}^{\infty}\frac{k!}{(k-n)!}(k-n)\delta_{k}\eta^kz^{k-n}\phi(w)^{k-n-1}\phi'(w)\right).
\end{align*}
By letting $w=0$ in the above expression, we obtain
\begin{equation}\label{Vasu-Him-Sub-P2-e-2.24}
	\phi(z)=b+\frac{n+3}{(n+1)^2}c \cdot \frac{\sum_{k=n+1}^{\infty}\frac{k!}{(k-n-1)!}\delta_{k}(\eta b)^{k-n-1}z^k}{\sum_{k=n}^{\infty}\frac{k!}{(k-n)!}\delta_{k}(\eta b)^{k-n}z^k}\,\,\, \mbox{for all}\,\, z\in \mathbb{D},
\end{equation}
where $c=\phi'(0)$.\\[1mm]

Conversely, let $\psi$ and $\phi$ be respectively of the form given by \eqref{Vasu-Him-Sub-e-2.23} and \eqref{Vasu-Him-Sub-P2-e-2.24} for some $a,b,c \in \mathbb{C}$. Then it is easy to see that for any $z\in \mathbb{D}$,
\begin{equation*}
	\frac{F_2(z)}{F_1(z)}=\frac{(n+1)^2}{n+3}\cdot \sum_{j=1}^{\infty}d_j(\eta b)^{j-1}z^j
\end{equation*}
with $d_1=1$ and $d_j\in \mathbb{R}^{+}$ for $j=2, 3, \cdots .$ Therefore, $\phi(z)$ can be written as 
\begin{equation}\label{Vasu-Him-Sub-P2-e-2.25}
	\phi(z)=b+c\sum_{j=1}^{\infty}d_j(\eta b)^{j-1}z^j.
\end{equation}

For any $n\in \mathbb{N}$, the operator $D_{n, \psi, \phi}$ will be complex symmetric with conjugation $C_{\mu, \eta}$ on $\mathcal{S}^2_{1}(\mathbb{D})$ if \eqref{Vasu-Him-Sub-P2-e-2.21} holds for all $z, w\in \mathbb{D}$.\\[3mm]

In view of \eqref{Vasu-Him-Sub-e-2.23} and \eqref{Vasu-Him-Sub-P2-e-2.25}, the relation \eqref{Vasu-Him-Sub-P2-e-2.21} is equivalent to
\begin{align*}
	&	\left(\frac{a}{n!\delta_{n}}\sum_{k=n}^{\infty}\frac{k!}{(k-n)!}\delta_{k}(\eta b)^{k-n}z^k\right)\left(\sum_{m=n}^{\infty}\frac{m!}{(m-n)!}\delta_{m}(\eta w)^m\left(b+c\sum_{j=1}^{\infty}d_j(\eta b)^{j-1}z^j\right)^{m-n}\right)\\ \\[2mm]& =
	\left(\frac{a}{n!\delta_{n}}\sum_{k=n}^{\infty}\frac{k!}{(k-n)!}\delta_{k}(\eta b)^{k-n}w^k\right)\left(\sum_{m=n}^{\infty}\frac{m!}{(m-n)!}\delta_{m}(\eta z)^m\left(b+c\sum_{j=1}^{\infty}d_j(\eta b)^{j-1}w^j\right)^{m-n}\right)
\end{align*}
which implies
\begin{align*}
	& \left(\sum_{k=n}^{\infty}\frac{k!}{(k-n)!}\delta_{k}(\eta b)^{k-n}z^k\right)\left\{\sum_{m=n}^{\infty}\frac{m!}{(m-n)!}\delta_{m}(\eta w)^m \left(\sum_{l=0}^{m-n}b^l\left(\sum_{j=1}^{\infty}cd_j(\eta b)^{j-1}z^j\right)^{m-n-l}\right)\right\}\\ \\[2mm]& =
	\left(\sum_{k=n}^{\infty}\frac{k!}{(k-n)!}\delta_{k}(\eta b)^{k-n}w^k\right)\left\{\sum_{m=n}^{\infty}\frac{m!}{(m-n)!}\delta_{m}(\eta z)^m \left(\sum_{l=0}^{m-n}b^l\left(\sum_{j=1}^{\infty}cd_j(\eta b)^{j-1}w^j\right)^{m-n-l}\right)\right\}.
\end{align*}
By comparing the coefficients of $z^{n+2}w^{n+1}$ from the both sides of the above equation we obtain 
\begin{align}\label{Vasu-Him-Sub-P2-e-2.26}
	(n+1)!\eta^{n+1}\delta_{n+1}&\left(\frac{1}{2}(n+2)!\eta^2 b^3\delta_{n+1}+\frac{n!\delta^2_{n}}{(n+1)\delta_{n+1}}d_2\eta bc+n!d_1\eta bc \delta_{n}\right)\\\nonumber
	\\[1mm] &=\frac{1}{2}(n+2)!\eta^{n+2}\delta_{n+2}\left(\frac{2n!d_1bc\delta^2_{n}}{\delta_{n+1}}+ (n+1)!\eta b^3 \delta_{n+1}\right)\nonumber.
\end{align}
Here we notice that the relation \eqref{Vasu-Him-Sub-P2-e-2.26} holds only if, either $b=0$ or $c=0$ or both are zero. Now we discuss about the following cases.

{\it Case I:} Let $b=0$ and $c\neq 0$. Then, $\psi(z)=az^n$ and $\phi(z)=cz$. Therefore, in this case, 
\begin{equation*}
	D_{n, \psi, \phi}C_{\mu, \eta}K_w(z)=a\mu \sum_{k=n}^{\infty}\frac{k!}{(k-n)!\delta_kc^{k-n}(\eta zw)^k}= C_{\mu, \eta}	D^*_{n, \psi, \phi}K_w(z).
\end{equation*}

{\it Case II:} Let $b\neq 0$ and $c=0$. Then we have 
$$ \phi(z)=b \,\,\,\, \mbox{and}\,\,\,\, \psi(z)=\frac{a}{n!\delta_{n}}\sum_{k=n}^{\infty}\frac{k!}{(k-n)!}\delta_{k}z^k(\eta b)^{k-n}.$$
An observation shows that
\begin{align*}
	D_{n, \psi, \phi}C_{\mu, \eta}K_w(z) &=\frac{a\mu}{n!\delta_{n}}\sum_{k=n}^{\infty}\frac{k!}{(k-n)!}\delta_{k}z^k(\eta b)^{k-n}\sum_{k=n}^{\infty}\frac{k!}{(k-n)!}\delta_{k}(\eta w)^kb^{k-n}\\ &=C_{\mu, \eta}	D^*_{n, \psi, \phi}K_w(z).
\end{align*}

{\it Case III:} Let $b=c=0$. Then $\phi(z)=0$ and $\psi(z)=az^n$. Therefore,
\begin{equation*}
	D_{n, \psi, \phi}C_{\mu, \eta}K_w(z)= a\mu n!\delta_{n}(\eta zw)^n=C_{\mu, \eta}	D^*_{n, \psi, \phi}K_w(z).
\end{equation*}

Therefore, we conclude that $D_{n, \psi, \phi}$ is $C_{\mu, \eta}$-symmetric on $\mathcal{S}^2_{1}(\mathbb{D})$ for each case. This completes the proof.
\end{pf}

\noindent\textbf{Acknowledgment:} 
The first named author is supported by SERB-CRG and the third named author is supported by DST-INSPIRE Fellowship (IF 190721),  New Delhi, India.

\end{document}